\documentclass[a4paper,11pt]{article}
\usepackage{amsmath,amsthm,amssymb,amsfonts,fullpage}

\newtheorem{thm}{Theorem}[section]

\newtheorem{lem}[thm]{Lemma}
\newtheorem{propn}[thm]{Proposition}
\newtheorem{rem}[thm]{Remark}

\begin{document}

\title{Slow movement of a random walk\\on the range of a random walk\\in the presence of an external field}

\author{David Croydon\footnote{Dept of Statistics,
University of Warwick, Coventry CV4 7AL, UK;
{d.a.croydon@warwick.ac.uk}.}}
\maketitle

\begin{abstract}
In this article, a localisation result is proved for the biased random walk on the range of a simple random walk in high dimensions ($d\geq 5$). This demonstrates that, unlike in the supercritical percolation setting, a slowdown effect occurs as soon a non-trivial bias is introduced. The proof applies a decomposition of the underlying simple random walk path at its cut-times to relate the associated biased random walk to a one-dimensional random walk in a random environment in Sinai's regime.
\end{abstract}

\section{Introduction}

In studying random walks in random environments, there is a particular focus at the moment on understanding the effect of an external field. Indeed, some quite remarkable results have been proved in this area. For instance, whereas adding a deterministic unidirectional bias to the random walk on the integer lattice $\mathbb{Z}^d$ results in ballistic escape, the same has been shown not to hold for supercritical percolation clusters. Instead, the random environment arising in the percolation model creates traps which become stronger as the bias is increased, so that when the bias is set above a certain critical value, the speed of the biased random walk is zero \cite{BGP, SAS}. This phenomenon, which has also been observed for the biased random walks on supercritical Galton-Watson trees \cite{BFGH, LPP} and a one-dimensional percolation model \cite{AFH}, is of physical significance, as it helps to explain how a particle could in some circumstances actually move more slowly when the strength of an external field, such as gravity, is greater \cite{BD}.

For percolation on the integer lattice close to criticality, physicists have identified two potential trapping mechanisms for the associated biased random walk: `trapping in branches' and `traps along the backbone' \cite{BD}. More concretely, in high dimensions the incipient infinite cluster for bond percolation on the integer lattice is believed to be formed of a single infinite path -- the backbone, to which a collection of `branches' or `dangling ends' is attached. If the dangling end is aligned with the bias, then the random walk will find it easy to enter this section of the graph, but very difficult to escape. Similarly, there will be sections of the backbone that flow with and sections that flow against the bias, and this will mean the random walk will prefer to spend time in certain locations along it.

Given that rigourous results for the incipient infinite cluster for critical bond percolation in $\mathbb{Z}^d$ are currently rather limited, exploring the biased random walks on it directly is likely to be difficult. Nonetheless, the above heuristics motivate a number of interesting, but more tractable research problems, one of which will be the focus of this article. In particular, to investigate the effect of `traps along the backbone', it makes sense to initially study how the presence of an external field affects a random walk on a random path. A natural choice for such a path is the one generated by a simple random walk on $\mathbb{Z}^d$, and it is for this reason that we pursue here a study of the biased random walk on this object.

To state our main result, we first need to formally define a biased random walk on the range of a random walk. Let $(S_n)_{n\in\mathbb{Z}}$ be a two-sided random walk on $\mathbb{Z}^d$, i.e. suppose that $(S_n)_{n\geq 0}$ and $(S_{-n})_{n\geq 0}$ are independent random walks on $\mathbb{Z}^d$ starting from 0 built on a probability space with probability measure $\mathbf{P}$. The range of this process is defined to be the random graph $\mathcal{G}=(V(\mathcal{G}),E(\mathcal{G}))$ with vertex set
\begin{equation}\label{vertex}
V(\mathcal{G}):=\left\{S_n:n\in \mathbb{Z}\right\},
\end{equation}
and edge set
\begin{equation}\label{edge}
E(\mathcal{G}):=\left\{\{S_n,S_{n+1}\}:n\in \mathbb{Z}\right\}.
\end{equation}
Now, fix a bias parameter $\beta\geq 1$, and to each edge $e=\{e_-,e_+\}\in E(\mathcal{G})$, assign a conductance
\begin{equation}\label{cond}
\mu_e:=\beta^{\max\{e_-^{(1)},e_+^{(1)}\}},
\end{equation}
where $e^{(1)}_\pm$ is the first coordinate of $e_\pm$. The biased random walk on $\mathcal{G}$ is then the time-homogenous Markov chain $X=((X_n)_{n\geq 0}, \mathbf{P}_{x}^{\mathcal{G}},x\in V(\mathcal{G}))$ on $V(\mathcal{G})$ with transition probabilities
\[P_\mathcal{G}(x,y):=\frac{\mu_{\{x,y\}}}
{\mu(\{x\})},\]
where $\mu$ is a measure on $V(\mathcal{G})$ defined by $\mu(\{x\}):=\sum_{e\in {E}(\mathcal{G}):x\in e}\mu_e$. A simple check of the detailed balance equations shows that $\mu$ is the invariant measure for $X$. Note that, if $\beta$ is strictly greater than $1$, then the biased random walk $X$ prefers to move in the first coordinate direction. If, on the other hand, there is no bias, i.e. $\beta=1$, then the preceding definition leads to the usual simple random walk on $\mathcal{G}$. Finally, as is the usual terminology for random walks in random environments, for $x\in V(\mathcal{G})$, we say that $\mathbf{P}_{x}^\mathcal{G}$ is the quenched law of $X$ started from $x$. Since 0 is always an element of $V(\mathcal{G})$, we can also define an annealed law $\mathbb{P}$ for the biased random walk on $\mathcal{G}$ started from 0 by setting
\begin{equation}\label{anlaw}
\mathbb{P}:=\int \mathbf{P}_0^\mathcal{G}(\cdot){\rm d}\mathbf{P}.
\end{equation}
Under this law, we can prove the following theorem, which shows that, unlike the supercritical percolation case, any non-trivial value of the bias leads to a slowdown effect.

{\thm\label{srwthm} Fix a bias parameter $\beta>1$ and $d\geq 5$. If ${X}=({X}_n)_{n\geq0}$ is the biased random walk on the range ${\mathcal{G}}$ of the two-sided simple random walk $S$ in $\mathbb{Z}^d$, then there exists an $S$-measurable random variable $L_n$ taking values in $\mathbb{R}^d$ such that
\[\mathbb{P}\left(\left|\frac{X_n}{\log n}-L_n\right|>\varepsilon \right)\rightarrow 0,\]
for any $\varepsilon>0$. Moreover, $(L_n)_{n\geq 1}$ converges in distribution under $\mathbf{P}$ to a random variable $L_\beta$ whose distribution can be characterised explicitly.}

{\rem \label{remo} The characterisation of $L_\beta$ that will be given in the proof of Theorem \ref{srwthm} readily yields that the distribution of $L_\beta\log \beta$ is independent of $\beta$. Thus, as the bias is increased, the biased random walk will be found closer to the origin.}
\bigskip

To show that the unbiased random walk $X$ on the graph $\mathcal{G}$ in dimensions $d\geq 5$ is diffusive, it was exploited in \cite{rwrrw} that the point process of cut-times of $S$, i.e. those times where the past and future paths do not intersect (of which there are infinite), is stationary. In particular, this observation allowed $\mathcal{G}$ to be decomposed at cut-points into a stationary chain of finite graphs, effectively reducing the problem into a one-dimensional one. (Note that the same techniques are no longer applicable when $d\leq 4$, as there are no longer an infinite number of cut-times for the two-sided random walk path.) This idea will again prove useful when proving Theorem \ref{srwthm}, with the difference being that now it must be taken into account how the bias affects each of the graphs in the chain. Since the orientations of the graphs in the chain are random, it turns out that the one-dimensional model it is relevant to compare to is a random walk in a random environment in the so-called Sinai regime. It is now well-known that, because of the large traps that arise, a random walk in a random environment in Sinai's regime escapes at a rate $(\log n)^2$ \cite{Sinai}. This will also be true for $X$ with respect to the graph distance, but taking into account that $S$ satisfies a diffusive scaling, we arrive at the $\log n$ scaling of the result.

The main difficulty in pursuing this line of reasoning is that the underlying simple random walk $S$ has loops, and so it is necessary to estimate how much time the biased random walk $X$ spends in these. If we start from a random path that is non-self intersecting, then there is not such a problem and, as long as the first coordinate of the random path still converges to a Brownian motion, verifying that a biased random walk exhibits a localisation phenomenon is much more straightforward. Thus, as a warm up to proving Theorem \ref{srwthm}, we start by considering biased random walks on non-self intersecting paths. As a particular example, we are able to prove the following annealed scaling limit for the biased random walk on the range of a two-sided loop-erased random walk in high dimensions (see Section \ref{brwsec} for precise definitions).

{\thm\label{lerwthm} Fix a bias parameter $\beta>1$ and $d\geq 5$. If $\tilde{X}=(\tilde{X}_n)_{n\geq0}$ is the biased random walk on the range $\tilde{\mathcal{G}}$ of the two-sided loop-erased random walk $\tilde{S}$ in $\mathbb{Z}^d$, then there exists an $\tilde{S}$-measurable random variable $\tilde{L}_n$ taking values in $\mathbb{R}^d$ such that
\[\mathbb{P}\left(\left|\frac{\tilde{X}_n}{\log n}-\tilde{L}_n\right|>\varepsilon \right)\rightarrow 0,\]
for any $\varepsilon>0$. Moreover, $(\tilde{L}_n)_{n\geq 1}$ converges in distribution under $\mathbf{P}$ to a random variable $\tilde{L}_\beta$ whose distribution can be characterised explicitly.}
\bigskip

This article contains only two further sections. In Section \ref{brwsec} we explain the relationship between the biased random walk on a random path and a random walk in a one-dimensional random environment, and prove Theorem \ref{lerwthm}. In Section \ref{srwsec}, we adapt the argument in order to prove Theorem \ref{srwthm}.

\section{Biased random walk on a self-avoiding random path}\label{brwsec}

The aim of this section is to describe how a biased random walk on a self-avoiding random path can be expressed as a random walk in a one-dimensional random environment. As we will demonstrate, this enables us to transfer results proved for the latter model to the former. To illustrate this, we will apply our techniques to the biased random walk on the range of the two-sided loop-erased random walk in dimensions $d\geq 5$.

We start by introducing some notation. Suppose that $S=(S_n)_{n\in \mathbb{Z}}$ is a random self-avoiding path in $\mathbb{R}^d$ with $S_0=0$ built on a probability space with probability measure $\mathbf{P}$. The range of this process $\mathcal{G}=(V(\mathcal{G}),E(\mathcal{G}))$ is defined analogously to (\ref{vertex}) and (\ref{edge}), so that, by the self-avoiding assumption, $\mathcal{G}$ is a bi-infinite path. Assign edge conductances as at (\ref{cond}), and let $X=((X_n)_{n\geq 0}, \mathbf{P}_{x}^{\mathcal{G}},x\in V(\mathcal{G}))$ be the associated biased random walk, i.e. the time-homogenous Markov chain  on $V(\mathcal{G})$ with transition probabilities
\[P_\mathcal{G}(S_n,S_{n\pm 1}):=\frac{c(\{S_n,S_{n\pm 1}\})}
{c(\{S_n,S_{n- 1}\})+c(\{S_n,S_{n+1}\})}.\]
As well as the quenched laws $\mathbf{P}_{x}^\mathcal{G}$, we can also define the annealed law for the process $X$ started from 0 by integrating out the underlying random path $S$, cf. (\ref{anlaw}).

Now let us define the particular random walk in a random environment of interest to us in this section. Firstly, the random environment $\omega$ will be represented by a random sequence $(\omega_n^{-},\omega_n^{+})_{n\in\mathbb{Z}}$ in $[0,1]^2$ such that
$\omega^-_n+\omega_n^+=1$, and will again be built on the probability space with probability measure $\mathbf{P}$. The random walk in the random environment will be the  time-homogenous Markov chain ${X}'=(({X}'_n)_{n\geq 0}, \mathbf{P}_{x}^{\omega},x\in \mathbb{Z})$ on $\mathbb{Z}$ with transition probabilities
\[P_\omega(n,n\pm1 )=\omega^\pm_n.\]
For $x\in \mathbb{Z}$, the law $\mathbf{P}_x^\omega$ is the quenched law of ${X}'$ started from $x$. Moreover,  we can define an annealed law for ${X}'$ started from 0 by integrating out the environment, similarly to \eqref{anlaw}. To connect this model with the biased random walk on the random path introduced above, we suppose that the transition probabilities  are defined by setting $\omega^\pm_n=P_\mathcal{G}(S_n,S_{n\pm1})$.
For this choice of random environment, it is immediate that, for any $x\in V(\mathcal{G})$, the law $\mathbf{P}^\omega_x\circ S^{-1}$, where $S^{-1}$ is the pre-image of the map $n\mapsto S_n$, is precisely the same as $\mathbf{P}^\mathcal{G}_x$. In other words, the quenched law of $S_{{X}'}$ is the same as that of $X$. A corresponding identity holds for the relevant annealed laws.

Importantly, it is also possible to connect the  first coordinate of the random path with the potential of the random walk in the random environment. To be more concrete, let $(S_n^{(1)})_{n\in\mathbb{Z}}$ be the first coordinate of $(S_n)_{n\in\mathbb{Z}}$, and $(\Delta_n)_{n\in\mathbb{Z}}$ be its increment process, i.e. \[\Delta_n:=S^{(1)}_n-S^{(1)}_{n-1}.\]
Then, if $\rho_n:=\omega_n^-/\omega_n^+$, where $\omega$ is defined as in the previous paragraph, an elementary calculation yields $\log \rho_n = -\log\beta(\Delta_{n+1}^+-\Delta_n^-)$, where $\Delta_n^+:=\max\{0,\Delta_n\}$ and  $\Delta_n^-:=-\min\{0,\Delta_n\}$. Consequently, the potential $(R_n)_{n\in\mathbb{Z}}$ of the random walk in a random environment, which is obtained by setting
\begin{equation}\label{potential}
R_n:=\left\{\begin{array}{cc}
                \sum_{i=1}^n\log\rho_i, & \mbox{if }n\geq 1,\\
                0, & \mbox{if }n=0, \\
                -\sum_{i=n+1}^0\log\rho_i, & \mbox{if }n\leq -1,
              \end{array}\right.
              \end{equation}
satisfies
\begin{equation}\label{potentialinc}
R_n=-\log\beta\left(S_n^{(1)}+\Delta_{n+1}^+-\Delta_1^+\right).
\end{equation}
Hence, if the individual increments are small, the first coordinate of $S$ very nearly gives a (negative) constant multiple of the potential of the random walk in the random environment.

The potential is of particular relevance when understanding the behaviour of the random walk in a random environment in the Sinai regime. In particular, by applying the fact that the potential converges to a Brownian motion, it is possible to describe where the large traps in the environment appear, and thus where the random walk prefers to spend time. Hence, at least when $S$ satisfies a scaling result that incorporates a functional invariance principle in the first coordinate (and the increments of $S^{(1)}$ are bounded), it is possible to use the relationship between $S^{(1)}$ and $R$ derived above to obtain the behaviour of the biased random walk on the random path.

{\propn \label{main} Fix a bias parameter $\beta>1$. Suppose that $S$ satisfies
\[\left(n^{-1/2}S_{\lfloor nt\rfloor}\right)_{t\in\mathbb{R}}\rightarrow(B_t)_{t\in\mathbb{R}}\]
in distribution, where $(B_t)_{t\in\mathbb{R}}$ is a continuous $\mathbb{R}^d$-valued process whose first coordinate $(B^{(1)}_t)_{t\in\mathbb{R}}$ is a non-trivial multiple of a standard two-sided one-dimensional Brownian motion. Moreover, suppose that the increment process $(\Delta_n)_{n\in\mathbb{Z}}$ satisfies $|\Delta_0|<C$, $\mathbf{P}$-a.s., for some deterministic constant $C$. It then holds that the biased random walk $X$ satisfies
\[\mathbb{P}\left(\left|\frac{X_n}{\log n}-L_n\right|>\varepsilon\right)\rightarrow 0,\]
for any $\varepsilon>0$, where $L_n$ is an $S$-measurable random variable that converges in distribution under $\mathbf{P}$ to a non-trivial random variable $L_\beta$ whose distribution can be characterised explicitly.}
\begin{proof} Recalling the identity at (\ref{potentialinc}), it is clear that the assumptions on $S^{(1)}$ imply the potential $R$ converges when rescaled to a Brownian motion. Hence, by applying the proof of \cite[Theorem 2.5.3]{Zeitouni} (and the following discussion), it is possible to demonstrate that the random walk in the random environment $X'$ satisfies
\[\mathbb{P}\left(\left|\frac{X'_n}{(\log n)^2}-b(n)\right|>\varepsilon\right)\rightarrow 0,\]
for any $\varepsilon>0$, where $b(n)$ is an $S$-measurable random variable that converges in distribution under $\mathbf{P}$ to a non-trivial random variable $b$ whose distribution can be characterised explicitly. (Note that, in the case when $(\Delta_n)_{n\in\mathbb{Z}}$ is an i.i.d. sequence, this result follows from \cite{Kestenrwre} and \cite{Sinai}.) Setting
\[L_n:=\frac{S_{(\log n)^2b(n)}}{\log n},\]
and $L:=B_b$, the proposition readily follows.
\end{proof}

To conclude this section, we note that the above result applies when $S$ is a two-sided loop-erased random walk in dimension $d\geq 5$. To introduce this model, we follow \cite[Chapter 7]{Lawlerbook}. First, fix $d\geq 5$ and suppose that $(\xi_n)_{n\geq 0}$ is a simple random walk on the integer lattice $\mathbb{Z}^d$. By the transience of this process, it is possible to define a sequence $(\sigma_n)_{n\geq 0}$ by setting $\sigma_0=0$ and, for $n\geq 1$, $\sigma_n:=\sup\{m: \xi_m=\xi_{\sigma_{n-1}+1}\}$. The loop-erasure of $(\xi_n)_{n\geq 0}$ is then the process $(S'_n)_{n\geq 0}$, where $S'_n:=\xi_{\sigma_n}$. Roughly speaking, $S'$ is derived from $\xi$ by erasing the loops of the latter process in a chronological order. To construct a two-sided version of the loop-erased random walk, we now suppose that we have two independent random walks on $\mathbb{Z}^d$ started from the origin, $\xi^1$ and $\xi^2$ say. Let $S^1$, $S^2$ be the loop-erasures of $\xi^1$, $\xi^2$, respectively, and set $A:=\{\xi^1_{[0,\infty)}\cap \xi^2_{[1,\infty)}=\emptyset\}$, which is an event with strictly positive probability in the dimensions that we are considering. On the event $A$, we then define $(S_{n})_{n\in\mathbb{Z}}$ by setting
\[S_n:=\left\{\begin{array}{cc}
                S^1_{-n}, & \mbox{if }n\leq 0,\\
                S^2_n, & \mbox{if }n\geq 0.
              \end{array}\right.\]
The process $S$ under the conditional law $\mathbf{P}(\cdot|A)$, where $\mathbf{P}$ is the probability measure on the space on which the two original random walks are defined, is the two-sided loop-erased random walk. Note that this is the same as the process defined in \cite{Lawlerlerw}, Section 5. Since $S$ is a nearest-neighbour path in $\mathbb{Z}^d$, the corresponding increments $(\Delta_n)_{n\in\mathbb{Z}}$ are clearly bounded. Moreover, that $(n^{-1/2}S_{\lfloor nt\rfloor})_{t\in\mathbb{R}}$ converges to a $d$-dimensional Brownian motion is effectively proved in \cite{Lawlerlerw}, Section 5. Thus the assumptions of Proposition \ref{main} are satisfied, and Theorem \ref{lerwthm} follows.

\section{Biased random walk on the range of simple random walk}\label{srwsec}

The goal of this section is to develop the techniques of the previous section to deduce results about the biased random walk on the range of a two-sided simple random walk in high dimensions. As noted in the introduction, the extra difficulty is that the underlying simple random walk has self-intersections, and so the range is no longer a simple path.

We start by introducing the notation that will allow us to study the biased random walk observed at the hitting times of cut-points. Let
\[\mathcal{T}:=\left\{n:S_{(-\infty,n]}\cap S_{[n+1,\infty)}=\emptyset\right\}\]
be the set of cut-times for $S$. This set is infinite, $\mathbf{P}$-a.s., and so we can write $\mathcal{T}=\{T_n:n\in\mathbb{Z}\}$, where $\dots <T_{-1}< T_0\leq0< T_1<T_2<\dots$. The corresponding cut-points will be denoted $C_n:=S_{T_n}$. Define the hitting times by $X$ of the set of cut-points $\mathcal{C}:=\{C_n:n\in\mathbb{Z}\}$ by setting
\[H_0:=\inf\{m\geq 0: X_m\in\mathcal{C}\},\]
and, for $n\geq 0$,
\[H_n:=\inf \{m>H_{n-1}:X_m\in\mathcal{C}\}.\]
Denoting by $\pi$ the bijection from $\mathbb{Z}$ to $\mathcal{C}$ that satisfies $\pi(n)=C_n$, we then let $(J_n)_{n\geq 0}$ be the $\mathbb{Z}$-valued process obtained by setting
\begin{equation}\label{jdef}
J_n:=\pi^{-1}\left(X_{H_n}\right).
\end{equation}

The parallel with Section \ref{brwsec} is that $J$ is a random walk in a random environment.
Note that, unlike in Section \ref{brwsec}, we allow the possibility that $J$ sits at a particular integer for multiple time-steps, and to capture this we will now write the environment as ($\omega_n^-,\omega_n^0,\omega_n^+)_{n\in\mathbb{Z}}$, where $\omega_n^\pm$ are defined to be the jump probabilities to $n\pm1$ from $n$, and $\omega_n^0$ is the probability of remaining at $n$. In particular, it is a simple calculation to check that
\[\omega_n^{\pm}:=\frac{1}{\mu(\{C_n\})R_{\rm eff}(C_n,C_{n\pm 1})},\]
where $R_{\rm eff}$ is the effective resistance operator on $V(\mathcal{G})$ corresponding to the given conductances (cf. \cite[(12)]{rwrrw}). As in the previous section, we will write $\rho_n:=\omega_n^{-}/\omega_n^+$ and define from this a potential $(R_n)_{n\in\mathbb{Z}}$ as at (\ref{potential}). Our first step is to show that this potential satisfies a functional invariance principle.

{\lem \label{potlemma} Fix a bias parameter $\beta>1$ and $d\geq 5$. The potential of the random environment $\omega$ satisfies
\[\left(n^{-1/2}R_{\lfloor nt \rfloor}\right)_{t\in\mathbb{R}}\rightarrow \left(\sigma B_t\right)_{t\in\mathbb{R}}\]
in distribution under $\mathbf{P}$, where $(B_t)_{t\geq 0}$ is a standard two-sided one-dimensional Brownian motion with $B_0=0$ and
\[\sigma^2:=\frac{(\log \beta)^2\mathbf{E}(T_1|0\in\mathcal{T})}{d}\in(0,\infty).\]}
\begin{proof} We will start by showing that, similarly to (\ref{potentialinc}), $(R_n)_{n\in\mathbb{Z}}$ is close to a constant multiple of the first coordinate of the cut-time process $(C_n^{(1)})_{n\in\mathbb{Z}}$. We can write
\[\log \rho_n=\log\left(\frac{\omega_n^-}{\omega_n^+}\right)=\log {R_{\rm eff}(C_n,C_{n+1})}-\log{R_{\rm eff}(C_n,C_{n-1})}.\]
Hence,
\[R_n=\log R_{\rm eff}(C_n,C_{n+1})-\log R_{\rm eff}(C_0,C_1).\]
Noting that the effective resistance between two vertices is always less than the graph distance between them in the graph when edges are weighted according to their individual resistances, it is possible to deduce that
\begin{equation}\label{renn1}
R_{\rm eff}(C_n,C_{n+1}) \leq \sum_{m=T_n}^{T_{n+1}-1}c(\{S_m,S_{m+1}\})^{-1}\leq \sum_{m=T_n}^{T_{n+1}-1}\beta^{-S_m^{(1)}}\leq (T_{n+1}-T_{n})\sup_{T_n\leq m\leq T_{n+1}-1}\beta^{-S_m^{(1)}}.
\end{equation}
Furthermore, since any path from $C_n$ to $C_{n+1}$ must contain the edge $\{S_{T_{n}},S_{T_{n+1}}\}$, it also holds that
\begin{equation}\label{renn2}
R_{\rm eff}(C_n,C_{n+1})\geq c(\{S_{T_{n}},S_{T_{n}+1}\})^{-1}= \beta^{-\max\{S_{T_{n}}^{(1)},S_{T_{n}+1}^{(1)}\}}.
\end{equation}
Thus,
\[\left|\log R_{\rm eff}(C_n,C_{n+1})+C_n^{(1)}\log \beta\right|\leq \log(T_{n+1}-T_n)+\log\beta\sup_{T_n\leq m\leq T_{n+1}}\left|C_n^{(1)}-S_m^{(1)}\right|,\]
and so
\begin{equation}\label{esti}
{\sup_{|m|\leq n}\left|R_m+C_m^{(1)}\log\beta\right|}\leq 2\sup_{|m|\leq n}\left[\log(T_{m+1}-T_m)+\log\beta\sup_{T_m\leq k\leq T_{m+1}}\left|C_m^{(1)}-S_k^{(1)}\right|\right].
\end{equation}

By a simple time-change, the estimate of the previous paragraph will allow us to prove the lemma from the obvious invariance principle for the first coordinate of the random walk,
\begin{equation}\label{sconv}
\left(n^{-1/2}S^{(1)}_{\lfloor nt\rfloor}\right)_{t\in\mathbb{R}}\rightarrow \left(d^{-1/2}B_{t}\right)_{t\in\mathbb{R}}.
\end{equation}
In particular, an ergodic theory argument implies that $n^{-1}T_n\rightarrow {\mathbf{E}}(T_1|0\in\mathcal{T})\in[1,\infty)$ as $|n|\rightarrow \infty$ almost-surely with respect to ${\mathbf{P}}(\cdot|0\in\mathcal{T})$ (see \cite[Lemma 2.2]{rwrrw}), and that the same holds true $\mathbf{P}$-a.s. can be shown by applying the relationship between the conditioned and unconditioned measures of \cite[(1.11)]{BSZ}. It readily follows that
\[\left(n^{-1/2}C^{(1)}_{\lfloor nt\rfloor}\log \beta\right)_{t\in\mathbb{R}}\rightarrow \left(\sigma B_t\right)_{t\in\mathbb{R}},\]
and so to complete the proof it will suffice to show that, when rescaled by $n^{-1/2}$, the right-hand side of (\ref{esti}) converges to 0 in $\mathbf{P}$-probability. To prove this, first observe that since $n^{-1}T_n$ converges $\mathbf{P}$-a.s., we further have
\[n^{-1}\sup_{|m|\leq n}(T_{m+1}-T_m)\rightarrow 0,\]
$\mathbf{P}$-a.s. The relevant convergence can be deduced from this and the tightness of $(n^{-1/2}S^{(1)}_{\lfloor nt\rfloor})_{t\in\mathbb{R}}$ that is an immediate consequence of (\ref{sconv}).
\end{proof}

To introduce the valleys of the potential, which play an important role in determining the behaviour of the random walk, we follow the presentation of \cite[Section 2.5]{Zeitouni}. A triple $(a,b,c)\in \mathbb{Z}^3$ with $a<b<c$ is a valley of $R$ if
\[R_a=\max_{a\leq n\leq b} R_n,\hspace{20pt}R_b =\min_{a\leq n\leq c} R_n,\hspace{20pt}R_{c} = \max_{b\leq n\leq c}R_n.\]
The depth of the valley is defined to be equal to
\[\min\left\{R_a-R_b,R_c-R_b\right\}.\]
If  $(a,b,c)$ is a valley of $R$ and $a<d<e<b$ are such that
\[R_e-R_d=\max_{a\leq m<n\leq b}(R_n-R_m),\]
then $(a,d,e)$ and $(e,b,c)$ are again valleys, obtained from $(a,b,c)$ by a so-called left-refinement. One can similarly define a right-refinement. Now, for $n\geq 2$, let
\begin{eqnarray*}
a'(n)&:=&\sup\{m \leq 0: R_m\geq \log n\},\\
c'(n)&:=&\inf\{m\geq 0: R_m\geq \log n\}
\end{eqnarray*}
and $b'(n)$ be the smallest integer in $[a'(n),c'(n)]$ where $R_{b'(n)}=\min_{a'(n)\leq m\leq c'(n)}R_m$, so that $(a'(n),b'(n),c'(n))$ is a valley of $R$ of depth $\geq \log n$.  By taking a successive sequence of refinements of $(a'(n),b'(n),c'(n))$, we can find the `smallest' valley $(a(n),b(n),c(n))$ with $a(n)<0$, $c(n)>0$ and depth $\geq \log n$. For $\delta>0$, the smallest valley $(a_\delta(n),b_\delta(n),c_\delta(n))$ with depth $\geq (1+\delta)\log n$ is defined similarly.

In much of what follows, it will be useful to assume that the random environment satisfies certain properties. To this end, we define $A(n,K,\delta)$ to be the subset of the probability space on which the random walk $S$ is built where:
\begin{itemize}
   \item $b(n)=b_\delta(n)$,
   \item any refinement $(a,b,c)$ of $(a_\delta(n),b_\delta(n),c_\delta(n))$ with $b \neq b(n)$ has depth $<(1-\delta)\log n$,
   \item $\min_{m\in[a_\delta(n),c_\delta(n)]\backslash [b(n)-\delta(\log n)^2,b(n)+\delta(\log n)^2]}(R_m-R_{b(n)})>\delta^3\log n$,
   \item $|a_\delta(n)|+|c_\delta(n)|\leq K(\log n)^2$,
   \item $\sup_{|m|\leq K(\log n)^2+1}\left[\log(T_{m+1}-T_m)+\log\beta\sup_{T_m\leq k\leq T_{m+1}}\left|C_m^{(1)}-S_k^{(1)}\right|\right]\leq {\delta^4}\log n$.
\end{itemize}
We note that
\[\lim_{\delta\rightarrow0}\limsup_{K\rightarrow\infty}\limsup_{n\rightarrow\infty}\mathbf{P}(A(n,K,\delta))=1.\]
Indeed, if we eliminate the final property, then this is essentially a restatement of  \cite[(2.5.2)]{Zeitouni}, and only depends on the fact that $R$ converges when rescaled to a Brownian motion. That we can incorporate the final property was verified in the proof of the previous lemma (with $n$ in place of $\log n$).

Before proceeding, we first observe that on $A(n,K,\delta)$ it is possible to derive a lower bound for the jump probabilities of the process $J$. More specifically, we claim that on the set in question
\begin{equation}\label{cond3}
\inf_{|m|\leq K(\log n)^2} \min\{\omega_m^-,\omega_m^+\}\geq (2d\beta)^{-1}n^{-2\delta^4}.
\end{equation}
To prove this, we apply the inequality at (\ref{renn1}) and the straightforward estimate $\mu(\{C_n\})\leq 2d\beta^{C_{n}^{(1)}+1}$ to obtain
\[\log\omega_n^{+}\geq -\log(2d\beta)-\log(T_{n+1}-T_n)-\log\beta \sup_{T_n\leq m\leq T_{n+1}}(C_n^{(1)}-S_m^{(1)}).\]
Since a similar lower bound also holds for $\log \omega_n^{-}$, the statement at (\ref{cond3}) follows from the final defining property of $A(n,K,\delta)$.

The following lemma outlines some first properties of the jump process $J$ defined at (\ref{jdef}).

{\lem \label{jlem} Fix a bias parameter $\beta>1$ and $d\geq 5$. For $\delta$ small and $K\in(0,\infty)$, there exists a finite integer $n_0(K,\delta)$ such that: if $n\geq n_0(K,\delta)$, then on $A(n,K,\delta)$ the jump process $J$ satisfies
\begin{equation}\label{firstest}
\mathbf{P}_0^\mathcal{G}\left(J\mbox{ hits $b(n)$ before time $\lfloor n^{1-\delta^2}\rfloor$}\right)\geq 1- n^{-\delta/4},
\end{equation}
and also
\begin{equation}\label{secondest}
\mathbf{P}_0^\mathcal{G}\left(\sup_{m\leq n}|J_m|\leq K(\log n)^2\right)\geq 1- n^{-\delta/4}.\end{equation}}
\begin{proof} For the first estimate, let us assume that $b(n)>0$. (The case $b(n)<0$ is similar, and the case $b(n)=0$ is trivial.) It is then a simple exercise in harmonic calculus to check that
\[\mathbf{P}_0^\mathcal{G}\left(J\mbox{ hits $a_\delta(n)$ before $b(n)$}\right)
\leq \frac{R_{\rm eff}(C_0,C_{b(n)})}{R_{\rm eff}(C_{a_\delta(n)},C_{b(n)})}=
\frac{\sum_{m=0}^{b(n)-1}R_{\rm eff}(C_m,C_{m+1})}{\sum_{m=a_\delta(n)}^{b(n)-1}R_{\rm eff}(C_{m},C_{m+1})},\]
where the inequality takes account of the fact that $J$ could start from $0$ or from $1$ if $X$ starts from 0. By applying the estimates for the effective resistance between cut-times from (\ref{renn1}) and (\ref{renn2}), the estimate for $R_m$ at (\ref{esti}) and the bounds that are known to hold on $A(n,K,\delta)$, it follows that
\begin{eqnarray*}
\lefteqn{\mathbf{P}_0^\mathcal{G}\left(J\mbox{ hits $a_\delta(n)$ before $b(n)$}\right)}\\
&\leq&{\frac{\sum_{m=0}^{b(n)-1}(T_{m+1}-T_m)\sup_{T_m\leq k\leq T_{m+1}}\beta^{C_m^{(1)}-S_k^{(1)}}e^{-R_m-C_m^{(1)}\log \beta}e^{R_m}}{\sum_{m=a_\delta(n)}^{b(n)-1}\beta^{-C_m^{(1)}-1}}}\\
&\leq&\beta  b(n)n^{5\delta^4}e^{\sup_{m\in[0,b(n)]}(R_m-R_{a_\delta(n)})}\\
&\leq&  \beta K (\log n)^2n^{5\delta^4-\delta}\\
&\leq&   n^{-\delta/2}
\end{eqnarray*}
for $\delta$ suitably small and $n\geq n_0(K,\delta)$. Furthermore, by proceeding as in the proof of \cite[Theorem 2.5.3]{Zeitouni}, it is possible to check that the expected time for the jump chain to hit the set $\{a_\delta(n),b(n)\}$ is bounded above by
\[\sum_{m=1}^{b(n)}\sum_{k=0}^{m-1-a_\delta(n)}\frac{R_{\rm eff}(C_{m-1},C_m) }{\omega_{m-k-1}^+R_{\rm eff}(C_{m-k},C_{m-k-1})}.\]
In turn, this can be bounded above by
\[4d\beta^2n^{9\delta^4} \sum_{m=1}^{b(n)}\sum_{k=0}^{m-1-a_\delta(n)}e^{R_m-R_{m-k}}\leq 4d \beta^2  K^2 (\log n)^4 n^{9\delta^4}e^{(1-\delta)\log n}\leq  n^{1-\delta/2},\]
again for $\delta$ chosen suitably small and $n\geq n_0(K,\delta)$, where, in addition to the estimates applied in the first part of the proof and the defining properties of $A(n,K,\delta)$, we have used the lower estimate for the transition probabilities from (\ref{cond3}). It is thus possible to conclude that
\begin{eqnarray*}
\lefteqn{\mathbf{P}_0^\mathcal{G}\left(J\mbox{ does not hit $b(n)$ before time $\lfloor n^{1-\delta^2}\rfloor$}\right)}\\
&\leq & \mathbf{P}_0^\mathcal{G}\left(J\mbox{ hits $a_\delta(n)$ before $b(n)$}\right)+ \mathbf{P}_0^\mathcal{G}\left(J\mbox{ does not hit $\{a_{\delta(n)},b(n)\}$ before time $\lfloor n^{1-\delta^2}\rfloor$}\right)\\
&\leq & n^{-\delta/2}+ \frac{n^{1-\delta/2}}{n^{1-\delta^2}}\\
&\leq & n^{-\delta/4}
\end{eqnarray*}
for small $\delta$ and $n\geq n_0(K,\delta)$, which completes the proof of (\ref{firstest}).

To prove (\ref{secondest}), we first observe that a similar argument to above yields
\[\mathbf{P}_{J_0=b(n)-1}^\mathcal{G}\left(J\mbox{ hits $a_\delta(n)$ before $b(n)$}\right)= \frac{R_{\rm eff}(C_{b(n)-1},C_{b(n)})}{R_{\rm eff}(C_{a_\delta(n)},C_{b(n)})}\leq \beta n^{6\delta^4}e^{R_{b(n)}-R_{a_\delta(n)}}\leq n^{-(1+\delta/2)}.\]
Similarly,
\[\mathbf{P}_{J_0=b(n)+1}^\mathcal{G}\left(J\mbox{ hits $c_\delta(n)$ before $b(n)$}\right)\leq n^{-(1+\delta/2)}.\]
Hence,
\begin{eqnarray}
\lefteqn{\mathbf{P}_{J_0=b(n)}^\mathcal{G}\left(\sup_{m\leq n}|J_m|\leq K(\log n)^2\right)}\nonumber\\
&\geq &\mathbf{P}_{J_0=b(n)}^\mathcal{G}\left(J\mbox{ returns to $b(n)$ at least $n$ times before hitting $\{a_\delta(n),c_\delta(n)\}$}\right)\nonumber\\
&\geq & \left(1-n^{-(1+\delta/2)}\right)^{n}\nonumber\\
& \geq & 1- n^{-\delta/3},\label{lower}
\end{eqnarray}
for $\delta$ small and $n\geq n_0(K,\delta)$. Since we also have that $J$ hits $b(n)$ before $\{a_\delta(n),c_\delta(n)\}$ with probability no less than $1-n^{\delta/2}$, the result follows.
\end{proof}

We now provide an upper estimate for the growth of hitting times.

{\lem \label{hlem} Fix a bias parameter $\beta>1$ and $d\geq 5$.  For $\delta$ small and $K\in(0,\infty)$, there exists a finite integer $n_0(K,\delta)$ such that: if $n\geq n_0(K,\delta)$, then on $A(n,K,\delta)$ the hitting time process $H$ satisfies
\[\mathbf{P}_0^\mathcal{G}\left( H_{\lfloor n^{1-\delta^2}\rfloor }\leq n\right)\geq 1- n^{-\delta^2/4}.\]}
\begin{proof} By simple properties of conditional expectation and the Markov property for $X$ (under the quenched law), we have that
\begin{eqnarray*}
\lefteqn{\mathbf{E}_0^\mathcal{G}\left(\left(H_{\lfloor n^{1-\delta^2}\rfloor}-H_0\right)\mathbf{1}_{\{\sup_{m\leq n}|J_m|\leq K(\log n)^2\}}\right)}\\
&= &\sum_{m=0}^{\lfloor n^{1-\delta^2}\rfloor-1}
\mathbf{E}_0^\mathcal{G}\left(\left(H_{m+1}-H_m\right)\mathbf{1}_{\{\sup_{m\leq n}|J_m|\leq K(\log n)^2\}}\right)\\
&\leq &\sum_{m=0}^{\lfloor n^{1-\delta^2}\rfloor-1}
\mathbf{E}_0^\mathcal{G}\left(\mathbf{E}_0^\mathcal{G}\left(H_{m+1}-H_m|\sigma(J_k:k\leq m)\right)\mathbf{1}_{\{|J_m|\leq K(\log n)^2\}}\right)\\
&=& \sum_{m=0}^{\lfloor n^{1-\delta^2}\rfloor-1}
\mathbf{E}_0^\mathcal{G}\left(\mathbf{E}_{C_{J_m}}^\mathcal{G}\left(H_{1}\right)\mathbf{1}_{\{|J_m|\leq K(\log n)^2\}}\right).
\end{eqnarray*}
Standard estimates for random walks on graphs in terms of volume and resistance (see \cite{Barlow}, Corollary 4.28, for example) imply  that the inner expectation satisfies
\begin{eqnarray*}
\mathbf{E}_{C_{J_m}}^\mathcal{G}\left(H_{1}\right)&\leq& R_{\rm eff}\left(C_{J_m},\{C_{J_m-1},C_{J_m+1}\}\right)\mu \left(\{S_k:T_{J_m-1}\leq k\leq T_{J_m+1}\}\right)\\
&\leq& R_{\rm eff}\left(C_{J_m},C_{J_m+1}\right)\sum_{k=T_{J_m-1}}^{T_{J_m+1}}2d\beta^{S_k^{(1)}+1}.
\end{eqnarray*}
Thus, on the set $\{|J_m|\leq K(\log n)^2\}$ it holds that
\[\mathbf{E}_{C_{J_m}}^\mathcal{G}\left(H_{1}\right)\leq 4d\beta n^{7\delta^4},\]
and so
\[{\mathbf{E}_0^\mathcal{G}\left(\left(H_{\lfloor n^{1-\delta^2}\rfloor}-H_0\right)\mathbf{1}_{\{\sup_{m\leq n}|J_m|\leq K(\log n)^2\}}\right)}\leq n^{1-\delta^2/2},\]
for small $\delta$ and $n\geq n_0(K,\delta)$. In fact, because one can similarly check that $\mathbf{E}_0^\mathcal{G}\left(H_0\right)\leq 2d\beta n^{6\delta^4}$, it is possible to replace $H_{\lfloor n^{1-\delta^2}\rfloor}-H_0$ by just $H_{\lfloor n^{1-\delta^2}\rfloor}$ in the above inequality. Consequently,
\[\mathbf{P}_0^\mathcal{G}\left( H_{\lfloor n^{1-\delta^2}\rfloor }> n,\sup_{m\leq n}|J_m|\leq K(\log n)^2\right)\leq n^{-\delta^2/2}.\]
In conjunction with (\ref{secondest}), this implies the result.
\end{proof}

All the pieces are now in place to prove Theorem \ref{srwthm} with
\[L_n:=\frac{C_{b(n)}}{\log n}.\]

\begin{proof}[Proof of Theorem \ref{srwthm}] As in the proof of \cite[Theorem 2.5.3]{Zeitouni}, the proof strategy will be to show that $X$ hits $C_{b(n)}$ before time $n$ and then stays there for a sufficient amount of time. For the majority of the proof, we will assume that $A(n,K,\delta)$ holds, with $\delta$ small and $n\geq n_0(K,\delta)$.

To show that $X$ hits $C_{b(n)}$ sufficiently early, we first observe that, by construction
\begin{eqnarray*}
\lefteqn{\mathbf{P}^\mathcal{G}_0\left(X\mbox{ doesn't hit $C_{b(n)}$ before time $n$}\right)}\\
&\leq&
\mathbf{P}^\mathcal{G}_0\left(J\mbox{ doesn't hit ${b(n)}$ before time $\lfloor n^{1-\delta^2}\rfloor$}\right) +
\mathbf{P}^\mathcal{G}_0\left(H_{\lfloor n^{1-\delta^2}\rfloor}> n\right).
\end{eqnarray*}
Hence, Lemmas \ref{jlem} and \ref{hlem} imply
\[\mathbf{P}^\mathcal{G}_0\left(X\mbox{ hits $C_{b(n)}$ before time $n$}\right)\geq 1-n^{-\delta^2/8},\]
for small $\delta$ and $n\geq n_0(K,\delta)$.

Now, since $J$ is the process $X$ observed at hitting times of the cut-point set $\mathcal{C}$, we are immediately able to deduce from (\ref{lower}) that
\[\mathbf{P}_{C_{b(n)}}^\mathcal{G}\left(X\mbox{ hits $\{C_{a_\delta(n)},C_{c_\delta(n)}\}$ before time $n$}\right)\leq n^{-\delta/3}.\]
It follows that
\[{\mathbf{P}_0^\mathcal{G}\left(\left|\frac{X_n}{\log n}-L_n\right|>\varepsilon\right)}\leq  n^{-\delta^2/8}+n^{-\delta/3}+\max_{m\leq n}\mathbf{P}_{C_{b(n)}}^\mathcal{G}\left(\left|\frac{\bar{X}_m}{\log n}-L_n\right|>\varepsilon\right),\]
where $\bar{X}$ is the random walk on the weighted graph $\bar{\mathcal{G}}$ with vertex set
\[V(\bar{\mathcal{G}}):=\{S_k:T_{a_\delta(n)}\leq k\leq T_{c_\delta(n)}\},\]
edge set
\[E(\bar{\mathcal{G}}):=\{\{S_k,S_{k+1}:T_{a_\delta(n)}\leq k\leq T_{c_\delta(n)}-1\},\] and edge conductances given by $\bar{c}(e)=c(e)$ (recall that $c(e)$ is the conductance of the edge $e$ in the original graph $\mathcal{G}$). To estimate the latter probability, we study the invariant measure $\bar{\mu}$ of $\bar{X}$, which is defined analogously to $\mu$. If $k\in[T_m,T_{m+1}]$, then
\[\bar{\mu}(\{S_k\})\leq 2d\beta^{S_k^{(1)}+1}\leq 2d\beta\sup_{k\in[T_m,T_{m+1}]}\beta^{S_k^{(1)}-C_m^{(1)}}e^{R_m+C_m^{(1)}\log\beta }e^{-R_m}.\]
Hence, if $m\in[a_\delta(n),c_\delta(n)]\backslash [b(n)-\delta(\log n)^2,b(n)+\delta(\log n)^2]$, then by applying (\ref{esti}) and the estimates that are known to hold on $A(n,K,\delta)$ it is possible to check that
\[\bar{\mu}(\{S_k\})\leq 2d\beta n^{3\delta^4-\delta^3}e^{-R_{b(n)}}.\]
Similarly, one can obtain
\[\bar{\mu}(\{C_{b(n)}\})\geq n^{-2\delta^4}e^{-R_{b(n)}}.\]
Since
\[\mathbf{1}_{C_{b(n)}}(x)\leq f(x):=\frac{\bar{\mu}(\{x\})}{\bar{\mu}(\{C_{b(n)}\})}, \hspace{20pt}\forall x\in V(\bar{\mathcal{G}}),\]
and $\bar{\mu}\bar{P}_\mathcal{G}=\bar{\mu}$, where $\bar{P}_\mathcal{G}$ is the transition matrix of $\bar{X}$, it follows that
\[\mathbf{P}_{C_{b(n)}}^\mathcal{G}\left(\bar{X}_l=S_k\right)=(\mathbf{1}_{C_{b(n)}}\bar{P}_\mathcal{G}^l)(S_k)\leq
(f\bar{P}_\mathcal{G}^l)(S_k)=f(S_k)\leq 2d\beta n^{5\delta^4-\delta^3}.\]
Thus,
\begin{eqnarray*}
\max_{m\leq n}\lefteqn{\mathbf{P}_{C_{b(n)}}^\mathcal{G}\left(\bar{X}_m = S_k\mbox{ for some }k\not\in[T_{b(n)-\delta(\log n)^2},T_{b(n)+\delta(\log n)^2}]\right)}\\
&\leq&
(T_{c_\delta(n)}-T_{a_\delta(n)})2d\beta n^{5\delta^4-\delta^3}\hspace{130pt}\\
&\leq &(c_\delta(n)-a_\delta(n))2d\beta n^{6\delta^4-\delta^3}\\
&\leq & 2d\beta K(\log n)^2 n^{6\delta^4-\delta^3}\\
&\leq &n^{-\delta^3/2}
\end{eqnarray*}
for small $\delta$ and $n\geq n_0(K,\delta)$.

We have thus reduced the problem to showing that
\begin{equation}\label{final}
\lim_{\delta\rightarrow0}\limsup_{n\rightarrow\infty}\mathbf{P}\left(\sup_{k\in[T_{b(n)-\delta(\log n)^2},T_{b(n)+\delta(\log n)^2}]}\left|\frac{S_k-C_{b(n)}}{\log n}\right|>\varepsilon\right)= 0.
\end{equation}
However, simultaneously with the convergence of $(d^{1/2}n^{-1/2}S_{\lfloor nt\rfloor})_{t\in\mathbb{R}}$ to a standard two-sided $d$-dimensional Brownian motion $(B_t)_{t\in\mathbb{R}}$ with $B_0=0$, one can check that $(\log n)^{-2}b(n)$ converges in distribution to some random variable $b(\infty)$ that takes values in $(-\infty,\infty)$ (cf. the discussion following \cite[Theorem 2.5.3]{Zeitouni}). Moreover, as was noted in the proof of Lemma \ref{potlemma}, $n^{-1}T_n$ converges $\mathbf{P}$-a.s. to a deterministic constant in $[1,\infty)$. Combining these results readily yields (\ref{final}).
\end{proof}

To complete the article, we will verify the statement of Remark \ref{remo} that $L_\beta\log \beta$, where $L_\beta$ is the distributional limit of $L_n$, has a distribution that is independent of $\beta$. Let $(B_t)_{t\in\mathbb{R}}$ be the standard two-sided Brownian motion that appears as the scaling limit of $(d^{1/2}n^{-1/2}S_{\lfloor nt\rfloor})_{t\in\mathbb{R}}$. Since $b(n)$ is the location of the base of the smallest valley of the process $(R_m)_{m\in\mathbb{Z}}$ that surrounds 0 and has depth $\log n$ , it is possible to check that $b(\infty)$, as defined in the previous proof, is the location of the base of the smallest valley of the process $(\frac{\log\beta}{\sqrt{d}}B^{(1)}_{t\tau})_{t\in\mathbb{R}}$, where $\tau:=\mathbf{E}(T_1|0\in\mathcal{T})$, which surrounds 0 and has depth 1. Moreover, $L_n$ converges to $L_\beta:=\frac{1}{\sqrt{d}}B_{b(\infty)\tau}$. By Brownian scaling, this implies that
\[L_\beta=\frac{B_{b'(\infty)}}{\log\beta}\]
in distribution, where $b'(\infty)$ is the location of the base of the smallest valley of $(B_t^{(1)})_{t\in\mathbb{R}}$ which surrounds 0 and has depth $1$, and so the claim does indeed hold true.

\def\cprime{$'$}
\providecommand{\bysame}{\leavevmode\hbox to3em{\hrulefill}\thinspace}
\providecommand{\MR}{\relax\ifhmode\unskip\space\fi MR }
\providecommand{\MRhref}[2]{%
  \href{http://www.ams.org/mathscinet-getitem?mr=#1}{#2}
}
\providecommand{\href}[2]{#2}

\end{document}